\documentclass[letterpaper,12pt]{amsart}
\usepackage[left=1in,right=1in,top=1in,bottom=1in]{geometry}

\usepackage{pgfplots}
\pgfplotsset{compat=1.18}

\usepackage{tikz, tikz-cd}
\usepackage{qtree}
\usepackage{forest}
\usepackage{enumerate}

\newcommand{\tinytree}[1]{
\vcenter{\hbox{
\scalebox{0.4}{
\begin{forest}
  for tree={circle,fill,}
  #1
\end{forest}
}}}}

\newcommand{\bigtree}[1]{
\vcenter{\hbox{
\scalebox{0.5}{
\begin{forest}
  for tree={circle,fill,s sep=5mm,inner sep=0.8mm}
  #1
\end{forest}
}}}}

\newcommand{\biggertree}[1]{
\vcenter{
\vspace{2mm}
\hbox{
\begin{forest}
  for tree={circle,fill,s sep=5mm,inner sep=0.8mm}
  #1
\end{forest}
\vspace{2mm}
}}}

\usepackage{graphicx}
\usepackage{float}
\usepackage{multicol}
\usepackage{soul}
\usepackage{pdfpages}

\tikzset{level distance=.7cm,sibling distance=0.5cm}
\usepackage{amsmath, amsfonts, amsthm, amssymb, mathtools}

\usetikzlibrary{shapes.geometric}

\theoremstyle{remark}
\newtheorem{example}{Example}[section]
\newtheorem{remark}{Remark}[section]

\theoremstyle{definition}
\newtheorem{definition}{Definition}[section]

\newtheorem{proposition}{Proposition}[section]

\newtheorem{theorem}{Theorem}[section]

\renewcommand{\tilde}{\widetilde}
\renewcommand*\bar[1]{\overline{#1}}

\renewcommand{\epsilon}{\varepsilon}

\renewcommand{\phi}{\varphi}

\newcommand{\eps}{\varepsilon}

\newcommand{\cA}{\mathcal{A}}
\newcommand{\cB}{\mathcal{B}}
\newcommand{\cC}{\mathcal{C}}
\newcommand{\cG}{\mathcal{G}}

\newcommand{\cH}{\mathcal{H}}
\newcommand{\cL}{\mathcal{L}}
\newcommand{\cM}{\mathcal{M}}

\newcommand{\cT}{\mathcal{T}}
\newcommand{\cU}{\mathcal{U}}

\DeclareFontFamily{U}{wncy}{}
\DeclareFontShape{U}{wncy}{m}{n}{<->wncyr10}{}
\DeclareSymbolFont{mcy}{U}{wncy}{m}{n}
\DeclareMathSymbol{\Sha}{\mathord}{mcy}{"58}

\newcommand{\K}{\mathbb{K}}

\newcommand{\fS}{\mathfrak{S}}

\newcommand{\fN}{\mathfrak{N}}
\newcommand{\fP}{\mathfrak{P}}

\newcommand{\fM}{\mathfrak{M}}

\DeclareMathOperator{\Prim}{Prim}

\DeclareMathOperator{\len}{len}
\DeclareMathOperator{\Sym}{Sym}
\DeclareMathOperator{\Hom}{Hom}
\DeclareMathOperator{\im}{im}

\newcommand{\ins}{\triangleleft}
\newcommand{\id}{\mathrm{id}}

\begin{document}

\title{Hopf algebras of leaf-graded forests}

\newcommand{\thisdate}{\today}

\author{Elizabeth Xiao}
\address{Department of Mathematics, California Institute of Technology, 1200 E California Blvd, Pasadena, CA 91125}
\email{exiao@caltech.edu}

\begin{abstract}
We examine the relationship between the Hopf algebra on nonplanar binary forests, previously defined by the author, and a Hopf algebra defined by Bruned and used in the renormalization of stochastic PDEs. We show that Bruned's Hopf algebra and its reduced analogue, previously studied by Foissy in the context of combinatorial twisted bialgebras, are dual to Hopf algebras derived from pre-Lie insertion operations respecting the number of leaves, which induce leaf-based gradings on the dual Hopf algebras. We compute Hochschild $1$-cocycles for these Hopf algebras and show that Bruned's Hopf algebra satisfies a similar universal property as the Connes-Kreimer Hopf algebra.
\end{abstract}

\maketitle

In \cite{xiao2026hopf} the author defined a Hopf algebra $\cH_{bin}$ on nonplanar rooted binary forests, which are the primary objects of study in Marcolli, Chomsky and Berwick's mathematical model of Merge in generative linguistics \cite{marcolli2025mathematical}. It was shown that the coproduct is indexed over binary-admissible cuts, a refinement of Connes-Kreimer's notion of admissible cut with the additional condition at most one of the edges below a vertex can be cut, and that the coproduct can also be indexed by subsets of leaves. It was also mentioned in Section 2.5 of \cite{xiao2026hopf} that $\cH_{bin}$ appears as the subalgebra of quotient of a renormalization Hopf algebra $\cH_B$ defined by Yvain Bruned in his PhD thesis \cite{bruned2015singular}, which is closely related to Connes-Kreimer's renormalization Hopf algebra $\cH_{CK}$.

Bruned's Hopf algebra $\cH_B$ is used in the renormalization of stochastic PDEs, and its precise combinatorial structures have not been fully studied until Foissy's recent work, first in \cite{foissy2021algebraic} and later in the monograph \cite{foissy2026combinatorial}, which places it in the context of combinatorial twisted bialgebras. We will show that both $\cH_B$ and its reduced analogue $\cH_{red}$ arise naturally from pre-Lie operators, much like $\cH_{bin}$. This provides another way to derive $H_{bin}$ as the restriction of a pre-Lie operator on $\cH_{red}$.

The paper is organized in the following manner. In Section 1, we fix notation, review the definition of a pre-Lie algebra, and remind the reader of the Guin-Oudom construction, which gives an explicit description of the universal enveloping algebra of a Lie algebra. In Section 2, we define what it means for a grading to be maximally compatible with a pre-Lie operation, which induces a restricted operation that still retains the Vinberg identity. In Section 3, we enumerate various pre-Lie insertion operations on trees. In Section 4, we show that $\cH_B$ is dual to the universal enveloping algebra obtained by restricting free pre-Lie insertion to a leaf-based grading; this follows directly from the definition of maximal compability. On the other hand, $\cH_{red}$ is related to a pre-Lie insertion on reduced trees which allows insertion at both edges and internal vertices; the duality of these algebras is shown in a combinatorial proof similar to the duality proofs in \cite{hoffman2003combinatorics} and \cite{xiao2026hopf}.

In Section 5, we review the definition of Hochschild cohomology for Hopf algebras. We show that $\cH_B$ has a 1-cocycle $B^\circ$ which is equal to the Connes-Kreimer 1-cocycle $B^+$, with the only difference being that $B^\circ(1)
 = 0$, instead of the single-leaf tree. We show that the pair $(\cH_B, B^\circ)$ satisfies a similar universal property as $(\cH_{CK}, B^+)$ with some additional conditions.

\section{Preliminaries}

\subsection{Rooted trees}

 A tree is a connected graph with no cycles. We assume that graphs are abstract in the sense that we do not fix a planar embedding. A nonplanar rooted tree is an abstract tree graph with a distinguished root vertex. We will draw all trees so that the root is at the top. Choosing a root induces a directed graph structure on a tree; all edges will be oriented to point away from the root.

A forest is a graph whose connected components are trees.

We will use the term tree to refer exclusively to nonplanar rooted trees; likewise, a forest will be a forest consisting of connected components which are nonplanar rooted trees. We can identify forests with multisets of trees, or commuting monomials of trees.

The degree or valency of a vertex is the number of the edges connected to it. We will often choose to use the term valency instead of degree when speaking of vertices, due to degree having another meaning in terms of graded vector spaces.

If $F$ is a forest, we write $V(F)$ for the vertex set of $F$; $E(F)$ for the edge set; $L(F) \subseteq V(F)$ for the subset of leaves, which, for trees containing more than one vertex, are non-root vertices of valency one. The root of a tree is never considered to be  a leaf if it has valency one. The only situation where the root is a leaf is the tree with a single vertex.

An internal vertex is any vertex which is not a leaf. We write $IV(F) \subseteq V(F)$ for the set of internal vertices. If the root of a tree has valency one, it is considered to be a internal vertex. (This is consistent with the definition, as the root in this case is not considered a leaf.)

Trees may have labels on all or some of its vertices and edges, although we do not consider labelled trees in this paper.

Given a vertex $v \in V(F)$, the vertex immediately above $v$ is called its parent, and the vertices below $v$ are called its children. Any vertex below $v$ is called a descendant of $v$. If $w \in V(F)$ is a child of $v$, we may refer to the whole subtree starting from $w$ as a child of $v$.

A tree is reduced if it contains no vertices of valency $2$, apart from possibly the root.

Let $\cT$ be some family of trees. Let $k$ be a field of characteristic $0$. We write $k\{\cT\}$ for the vector space over $k$ generated by the basis $\cT$. We identify $k[\cT]$ with $\Sym(k\{\cT\})$, the $k$-algebra of polynomials of trees, which is generated as a vector space by forests whose components are in $\cT$.

\subsection{Hopf algebras}
We refer to \cite{cartier2021classical} for a thorough description of algebras, coalgebras and Hopf algebras; this section summarizes the main terminology. Let $k$ be a field of characteristic $0$. A \emph{bialgebra} $\cB = (\cB, m, \eta, \Delta, \epsilon)$ over $k$ consists of the following data:
\begin{itemize}
  \item the underlying $k$-vector space $\cB$;
  \item a product $m : \cB \otimes \cB \to \cB$ with unit $\eta : k \to \cB$ (an unitary algebra structure);
  \item a coproduct $\Delta : \cB \to \cB \otimes \cB$ with counit $\epsilon : \cB \to k$ (a counitary coalgebra structure),
\end{itemize}
such that $\Delta$ and $\epsilon$ are morphisms of algebras. We often use \emph{Sweedler notation} to index the coproduct of an element $x \in \cB$ as a sum:
\[ \Delta(x) = \sum_{(x)} x' \otimes x''. \]

Since the unit and counit are assumed to exist, we can write $\cB$ as a triple $\cB = (\cB, m, \Delta)$. We also typically assume that $m$ is associative and $\Delta$ is coassociative.

We say $\cB$ is graded as a bialgebra if its underlying vector space is graded; $(\cB, m, \eta)$ is graded as an algebra; and $(\cB,\Delta,\eps)$ is graded as a coalgebra, in the sense that:
\begin{itemize}
  \item if $x \in \cB_n$ and $y \in \cB_m$, then $m(x,y) \in \cB_{n+m}$;
  \item if $x \in \cB_n$, then $\Delta(x) \in \bigoplus_{k=0}^n \cB_k \otimes \cB_{n-k}$.
\end{itemize}
A graded bialgebra $\cB$ is \emph{connected} if $\cB_0 \cong k$. It is \emph{locally finite} if $\dim \cB_n$ is finite for all $n$.

Given coalgebras $(\cC, \Delta, \eps)$ and $(\cC', \Delta', \eps')$, a coalgebra morphism is a map $\phi : \cC \to \cC'$ which satisfies
\[ (\phi \otimes \phi) \circ \Delta = \Delta ' \circ \phi \]
and
\[ \eps' \circ \phi = \eps. \]
A bialgebra morphism is a map which is both an algebra morphism and a coalgebra morphism.

We also remind the reader that the counit satisfies
\[ (\id \otimes \eps) \circ \Delta = \id = (\eps \otimes \id) \circ \Delta.  \]

A \emph{Hopf algebra} $\cH$ is a bialgebra with an antipode map $S : \cH \to \cH$ such that
\[ \eta(\epsilon(x)) =\sum_{(x)} S(x') x'' = \sum_{(x)} x' S(x'') \]
for all $x \in \cH$. If a bialgebra has an antipode, then it is unique.

A Hopf algebra morphism is any morphism between the underlying bialgebras; if $\cH, \cH'$ are Hopf algebras with antipodes $S, S'$, then any bialgebra map $\phi : \cH \to \cH'$ immediately satisfies $\phi \circ S = S' \circ \phi$.

If $\cB$ is a graded bialgebra, then it has an antipode defined recursively by $S(1) = 1$ and
\[ S(x) = -x - \sum_{(x)} S(x') x'' \]
for $x \in \bigoplus_{n=1}^\infty \cB_n$. Thus, any graded bialgebra is automatically a graded Hopf algebra.

Practically all Hopf algebras we work with are graded with counit $\eps : \cH \to k$ defined by $\eps(1) = 1$ and $\eps(h) = 0$ for $h \in \bigoplus_{n \geq 1} \cH_n$, which allows us to bypass verifying that Hopf algebra morphisms satisfy the counit condition.

\subsection{Pre-Lie algebras}

Let $\cL$ be a vector space. A bilinear product $\ins : \cL \times \cL \to \cL$ is called a \emph{(right) pre-Lie operator} if it satisfies the \emph{(right) Vinberg identity}
\[
   (x \ins y) \ins z - x \ins (y \ins z) =  (x \ins z) \ins y - x \ins (z \ins y),
\]
making $(\cL,\ins)$ a \emph{(right) pre-Lie algebra}.

The expression $(x \ins y) \ins z - x \ins (y \ins z)$ is  the \emph{associator} $A(x,y,z)$ with respect to $\ins$; if $\ins$ is associative, then $A(x,y,z) = 0$ uniformly. We write $A_\ins(x,y,z)$ when there is ambiguity in the operation of interest.

The Vinberg identity asserts that the associator is unchanged when the second and third arguments are swapped; that is, for all $x,y,z \in \cL$,
\[ A(x,y,z) = A(x,z,y). \]

A pre-Lie operator induces a Lie bracket defined by
\[ [x,y] := x \ins y - y \ins x, \]
which automatically satisfies antisymmetry and the Jacobi identity, making $\cL$ a Lie algebra. Not all Lie algebras, however, are produced from a pre-Lie algebra. Chapter 6 of \cite{cartier2021classical} provides a thorough exposition of pre-Lie algebras.

\subsection{The Guin-Oudom construction}

Guin and Oudom \cite{oudom2008lie} provide an explicit description of the universal enveloping algebra of a Lie algebra constructed from a pre-Lie operation, which is a Hopf algebra by the Poincar\'e-Birkhoff-Witt theorem. Another exposition of the Guin-Oudom procedure is found in \cite{foissy2025operads}.

Let $(\cL,\ins)$ be a pre-Lie algebra and let $\Sym(\cL)$ be the symmetric algebra of commuting polynomials over $\cL$. Then $\Sym(\cL)$ is naturally equipped with the shuffle coproduct, which is defined by
\[ \Sha(x_1 x_2 \dots x_n) = \sum_{I \subseteq [n]} x_I \otimes x_{[n] \setminus I} \]
where each $x_i \in \cL$ and, given an index set $J = \{j_1, \dots, j_k\}$, we write $x_J = x_{j_1} \cdots x_{j_k}$. The shuffle coproduct is coassociative and cocommutative.

The pre-Lie operation $\ins$ extends to an operation on $\Sym(\cL)$, which satisfies the following properties

\begin{enumerate}
  \item $A \ins 1 = A$
  \item $X \ins A Y = (X \ins A) \ins Y - X \ins (A \ins Y)$
  \item $AB \ins C = \sum_{(C)} (A \ins C_{(1)}) (B \ins C_{(2)})$
\end{enumerate}
for elements $X, Y \in \cL$ and $A,B,C \in \Sym(\cL)$. Recursively, this uniquely defines $A \ins B$ for any two monomials $A,B \in \Sym(\cL)$ and extends to all of $\Sym(\cL)$ by linearity.

The product of $\cU(\cL)$ as a Hopf algebra is a further modification of $\ins$. We form another product $\ast = \ast_\ins$ on $\Sym(\cL)$ by first taking the extension of $\ins$ to $\Sym(\cL)$ and then defining 
\[ A \ast B := \sum_{(B)} (A \ins B_{(1)} ) B_{(2)} \]
or monomials $A,B$, where the sum is indexed over terms of the shuffle coproduct (using Sweedler notation). Guin and Oudom prove that $\ast$ is associative.

This equips $(\Sym(\cL), \ast, \Sha)$ with the structure of a noncommutative and cocommutative Hopf algebra, which is isomorphic to the universal enveloping algebra of $\cL$. We will refer to this Hopf algebra as $\cU(\cL, \ins) = (\cU(\cL), \ast_\ins, \Sha)$, or even $\cU(\cL)$ when the pre-Lie operation is clear, as this is the only presentation of the universal enveloping algebra that we will explicitly use.

\section{Pre-Lie operations restricted by grading}

Given a pre-Lie algebra $(\cL, \ins)$, it is sometimes advantageous to consider different choices of grading on $\cL$ as a vector space. Let $\cL = \bigoplus_{n=0}^\infty \cL_n$ for some choice of grading, and let $p_k :\cL \to \cL_k$ be the projection to the $k$th graded component of $\cL$. We define a new operation $\ins'$ with respect to this grading by
\[ x \ins' y = p_{a+b}(x \ins y) \]
for homogeneous elements $x \in \cL_a$ and $y \in \cL_b$.

A pre-Lie algebra $(\cL, \ins)$ is graded if $\cL$ is graded and $x \ins y \in \cL_{a+b}$ for $x \in \cL_a$ and $y \in \cL_b$. If the operation $\ins'$ induced by grading is pre-Lie, then $(\cL, \ins')$ is graded as a pre-Lie algebra.

We say that a grading is \emph{maximally compatible} with $\ins$ if all the terms removed by $p_{a+b}$ are of lower degree than $a+b$, that is,
\[ x \ins y - x \ins'y \in \bigoplus_{k \leq a+b-1} \cL_k \]
for $x \in \cL_a$ and $y \in \cL_b$,
and \emph{minimally compatible} if all terms removed by $p_{a+b}$ are of higher degree than $a+b$. Note that maximal and minimal compatibility are essentially the same criterion if we allow negative gradings, so all proofs will be presented from the viewpoint of maximal compatibility.

\begin{remark}
  Not every grading of $\cL$ is compatible with $\ins$. For instance, one can assign an artificial grading to individual basis elements of $\cL$ and then cull terms from $x \ins' y$ which do not respect this grading, but there is no guarantee that $\ins'$ will still be pre-Lie.
\end{remark}

However, a finer notion called maximal (or minimal) compatibility does imply that the induced operator is pre-Lie.

\begin{proposition}
  Let $(\cL, \ins)$ be a pre-Lie algebra with a maximally compatible grading. Then the induced operation $\ins'$ is automatically pre-Lie.
\end{proposition}

\begin{proof}
  Suppose $L = \bigoplus_{n=0}^\infty \cL_n$ is a grading which is maximally compatible with $\ins$. Let $x \in \cL_a$, $y \in \cL_b$ and $z \in \cL_c$. Since $\ins$ is pre-Lie, it satisfies
  \[
    A_\ins(x,y,z)
    =
    A_\ins(x,z,y),
  \]
  and in particular, the degree $a+b+c$ components of the associators are equal:
  \[
    p_{a+b+c}(A_\ins(x,y,z))
    =
    p_{a+b+c}(A_\ins(x,z,y)).
  \]
  Consider now the associators of the induced operation. Since $\ins$ is maximally compatible with the grading, the terms of $x \ins y$, for instance, have degree no greater than $a+b$, and $x \ins' y$ keeps only the degree $a+b$ terms. Hence $A_\ins'(x,y,z)$ keeps only the degree $a+b+c$ terms of $A\ins(x,y,z)$, so
  \[ A_{\ins'}(x,y,z) = p_{a+b+c}(A_\ins(x,y,z)) \]
  and likewise for $A_{\ins'}(x,z,y)$, so
    \[
    A_{\ins'}(x,y,z)
    =
    A_{\ins'}(x,z,y),
  \]
  hence $\ins'$ is also pre-Lie.
\end{proof}

\subsection{The induced graded universal enveloping algebra}
Now we examine the universal enveloping algebra of a pair $(\cL, \ins')$ induced by a maximally compatible grading.

\begin{theorem}\label{compat}
Let $(\cL,\ins)$ be a pre-Lie algebra with a strongly compatible grading that induces the pair $(\cL,\ins')$. Suppose that $\cU(\cL, \ins)$ is dual to some Hopf algebra $(\cH,\Delta)$ with respect to a pairing $\langle -, - \rangle$. Then:

\begin{enumerate}
  \item there is an induced Hopf algebra $(\cH, \Delta')$ dual to $\cU(\cL, \ins')$, where the coproduct $\Delta'$ is a restriction of $\Delta$;

  \item there is a grading on $(\cH, \Delta')$, compatible with the coproduct, induced by the pre-Lie algebra;

  \item if $x \in \cH_n$, where $\cH = \bigoplus_{n=0}^\infty \cH_n$ is the grading induced in (b), then
\[ \Delta'(x) = \sum_{a+b = n} (( p_a \otimes p_b) \circ \Delta )(x). \]
\end{enumerate}
\end{theorem}

\begin{proof}
  A grading on $\cL$ extends to a grading on $\Sym(\cL)$ where the degree of a product $ X_1 \cdots X_n$ of homogeneous elements is the sum of the degrees of the individual factors,
  \[ \deg(X_1 \cdots X_n) = \sum_{i=1}^n \deg(X_i). \]
  In particular, concatenation of monomials respects degree. A monomial $X_1 \cdots X_n$ has length $n$; we denote this by $\len(X_1 \cdots X_n) = n$.

  We will show that the extension of $\ins$ to $\ast$ to monomials also respects degree, that is if $A, B \in \Sym(\cL)$ are two monomials of homogeneous elements of $\cL$, then
  \[ \deg(A \ast B) = \deg(A) + \deg(B). \]
  We first show that the extension of $\ins$ to monomials respects degree. First of all, recall that
  \begin{enumerate}[(i)]
  \item $A \ins 1 = A$
    \item $A \ins BX = (A \ins B) \ins X - A \ins (B \ins X)$
    \item $AB \ins C = \sum_{(C)} (A \ins C_{(1)}) (B \ins C_{(2)})$
  \end{enumerate}
  for $X, Y \in \cL$ and $A,B,C \in \Sym(\cL)$.
  
  We assume that $\deg(A' \ins B') = \deg A' + \deg B'$ holds whenever $A'$ is a product of at most $m$ homogeneous elements in $\cL$, and $B'$ is a product of at most $n$ homogeneous elements in $\cL$. Certainly this is true for $m=n=1$, since the original pre-Lie operation on $\cL$ respects grading.

  Note that we are inducting the \emph{length} of a monomial to draw conclusions about its \emph{degree}. Length and degree are both preserved by concatenation; however, length is not preserved by the extension of $\ins$. In fact, $\len(A \ins X) = \len A $.

  \begin{itemize}
   \item We will induct on the first argument of $\ins$. Suppose that $A$ is homogeneous of length $m$, while  while $B$ is homogeneous of length $n$. Let $Y \in \cL$ be homogeneous of length $1$. Then by the second property of the recursive extension,
   \[ AY \ins B = \sum_{(B)} (A \ins B_{(1)}) (Y \ins B_{(2)}), \]
   where $\len Y, \len A \leq m$, and for each pair $(B_{(1)}, B_{(2)})$,
   \[ \len B_{(i)} \leq  \len B \leq n \]
   for $i=1,2$. So the inductive hypothesis can be applied to obtain
   \[ \deg(A \ins B_{(1)}) = \deg A + \deg B_{(1)}  \]
   and
   \[ \deg(Y \ins B_{(2)}) = \deg Y + \deg B_{(2)}.  \]
   Since degree commutes with concatenation of polynomials, each term in $AY \ins B$ is homogeneous of degree
   \begin{align*}
     \deg((A \ins B_{(1)})(Y \ins B_{(2)}))
     &= (\deg A + \deg B_{(1)}) + (\deg Y + \deg B_{(2)}) \\
     &= \deg A + \deg Y + \deg B \\
     &= \deg (AY) + \deg B.
   \end{align*}
   Since $AY \ins B$ is a sum of homogeneous elements of the same degree,
   \[ \deg(AY \ins B) = \deg(AY) + \deg B, \]
   as required.

   \item We assume that $A$ is homogeneous of length $m$, while $B$ is homogeneous of length $n$. Let $X \in \cL$ be homogeneous of length $1$. Then
   \[ \deg(A \ins B) = \deg A + \deg B \]
   and the length of any term in $A \ins B$ is $\len A \leq m$ by property (3). Hence,
   \[ \deg((A \ins B) \ins X) = \deg(A \ins B) + \deg X = \deg A + \deg B + \deg X, \]
   since $\len(A \ins B) \leq m$ and $\len X = 1 \leq n$.

   From the induction in the first variable (even if $\deg B > m$) we have that
   \[ \deg(B \ins X) = \deg B + \deg X \]
   where $\len(B \ins X) = \len B \leq n$. Recall that $\deg A \leq m$ so the inductive hypothesis holds and
   \[ \deg(A \ins (B \ins X)) = \deg A + \deg (B \ins X) = \deg A + \deg B + \deg X \]
   Hence both terms in $A \ins BX$ are homogeneous of the same degree, so
   \[ \deg(A \ins BX) = \deg A + \deg (BX). \]

The extension of $\ins$ to $\Sym(\cL)$ is then modified to $\ast$ by defining
\[ A \ast B := \sum_{(B)} (A \ins B_{(1)}) B_{(2)} \]
where the degree of each term is
\begin{align*}
  \deg( (A \ins B_{(1)}) B_{(2)} )
  &= \deg(A \ins B_{(1)}) + \deg B_{(2)} \\
  &= \deg A + \deg B_{(1)} + \deg B_{(2)} \\
  &= \deg A + \deg B,
\end{align*}
so $\deg(A \ast B) = \deg A + \deg B$.
\end{itemize}

Suppose now that $(\cL, \ins)$ has a maximally compatible grading which induces the pair $(\cL, \ins')$. Let $\ast$ and $\ast'$ be the extensions, respectively, of $\ins$ and $\ins'$ to $\Sym(\cL)$. Since the extensions respect the grading, for homogeneous $A \in \Sym(\cL)_a$ and $B \in \Sym(\cL)_b$, the terms in $A \ast' B$ are exactly those in $A \ast B$ which have maximal degree $a+b$. That is,
\[ A \ast' B = p_{a+b}(A \ast B). \]

Now we are ready to prove (3), from which (1) and (2) immediately follow.

Write $(\cU(\cL), \ast, \Sha)$ as the dual of $(\cH, \sqcup, \Delta)$ with respect to the pairing $\langle -, - \rangle$. The proof of (3) is based on the duality between $\ast$ and $\Delta$, which is encapsulated by the following equation:
\[
  \langle A \ast B, C \rangle = \langle A \otimes B, \Delta(C) \rangle.
\]

Then $\Delta'$ is defined by
\[
  \langle A \ast' B, C \rangle = \langle A \otimes B, \Delta'(C) \rangle.
\]
Fixing $C \in \cH_n$ to be homogeneous, we observe that if $\langle A \ast B, C\rangle \neq 0$, then
\[ \langle A \ast' B, C \rangle = \langle A \ast B, C \rangle \]
exactly when $\deg A + \deg B = n$, and is equal to $0$ otherwise. Thus, if we write $\Delta(C) = \sum_{C', C''} C' \otimes C''$ in Sweedler notation, then
\begin{align*}
  \Delta'(C) = \sum_{
    \substack{
      C', C'' \\
      \deg C' + \deg C'' = n
      }
    } C' \otimes C''
    = \sum_{a+b = n} ( (p_a \otimes p_b) \circ \Delta)(C).
\end{align*}


\end{proof}

\section{Types of pre-Lie insertion}
Let $\cT$ be the collection of space of unlabelled rooted nonplanar trees. Let $\tilde \cT$ be the collection of reduced trees. Let $k\{\cT\}$ and $k\{\tilde \cT\}$ (resp. $k[\cT]$ and $k[\tilde \cT]$) be the vector spaces (resp. polynomial algebras) generated by each of these collections. Our pre-Lie operations will initially be defined on $k\{\cT\}$ and $k\{\tilde \cT\}$ and extended to the polynomial algebras through the Guin-Oudom procedure.

If $\cL$ is a vector space and $\ins$ is a pre-Lie operation on $\cL$, then we will write $\cU(\cL, \ins)$ for the universal enveloping algebra induced by this pair.

\subsection{The free pre-Lie algebra of vertex-insertion} The underlying vector space is $k\{ \cT\}$. Let $T, S \in \cT$. The operation $T \ins_{free} S$ enumerates all the trees obtained by choosing any vertex $v$ of $T$, leaves included, and joining $S$ to $v$ by creating a new edge from $v$ to the root of $S$.

\begin{example}
  Let $T = \tinytree{[[][[]]]}$ and $S = \bullet$. Then
  \[ T \ins S = \tinytree{[[][[]][]]} + \tinytree{[[[]][[]]]} + \tinytree{[[][[][]]]} + \tinytree{[[][[[]]]]}. \]
\end{example}

The operation $\ins_{free}$ respects grading with respect to vertex count. It was shown by Chapoton and Livernet in \cite{chapoton2001pre} that $\ins_{free}$ is a free pre-Lie operation on one generator, justifying the name we have chosen. 

The universal enveloping algebra $\cU(k\{\cT\}, \ins_{free})$ is isomorphic to the Grossman-Larson Hopf algebra, which was shown by Hoffman in \cite{hoffman2003combinatorics} to be dual to the Connes-Kreimer Hopf algebra. We delay further discussions of Hopf algebras and duality to later sections.

\subsection{Vertex-insertion with preservation of leaf count}

Let $T \in \cT$ be a tree with $n$ leaves. 
We introduce a new grading on $\cT$ such that trees are graded by the number of leaves. We assume that the root of $T \in \cT$ is not considered a leaf, even if it has a valency of $1$, unless $T = \bullet$. By this grading,
\begin{align*}
  \cT_1 &= \left\{
    \bullet,
    \tinytree{[[]]},
    \tinytree{[[[]]]},
    \dots
  \right\} \\
  \cT_2 &= \left\{
    \tinytree{[[][]]},
    \tinytree{[[[][]]]},
    \tinytree{[[[]][]]},
    \dots
  \right\} \\
  \cT_3 &= \left\{
    \tinytree{[[][][]]},
    \tinytree{[[[][]][]]},
    \tinytree{[[[][][]]]},
    \tinytree{[[[]][][]]},
    \dots
  \right\}
\end{align*}
This grading is not locally finite; every $\cT_n$ has infinite cardinality.

We define a new insertion operation $\ins_\ell$ which is the restriction of $\ins_{free}$ that keeps only the terms that respect the leaf-based grading. In practice, this simply removes all the terms in $T \ins_{free} S$ which are the result of inserting $S$ into a leaf of $T$, which causes the total number of leaves to decrease, meaning the leaf-grading is maximally compatible with $\ins_\ell$.

\begin{proposition}
  $\ins_\ell$ is a pre-Lie operation.
\end{proposition}

\begin{proof}
  The leaf-based grading is maximally compatible with $\ins_{free}$, inducing the pre-Lie operation $\ins_\ell$.  
\end{proof}

We are careful not to describe $\ins_\ell$ as ``insertion at internal vertices'', because the root of a tree is never considered an internal vertex, yet it is a possible insertion location for trees which are not a single leaf. This distinction will become irrelevant once we pass to reduced trees.

\begin{example}
  Let $T = \tinytree{[[][[]]]}$ and $S = \bullet$. Then
  \[ T \ins_\ell S = \tinytree{[[][[]][]]} + \tinytree{[[][[][]]]}. \]
\end{example}

\begin{remark}
The insertion $\ins_\ell$ is not free because $\bullet \ins_\ell T = 0$ for any tree $T$, since there are no valid insertion spots.
\end{remark}

\subsection{Edge-insertion on reduced trees}
Our underlying vector space is now $k\{\tilde \cT \}$. We will define $\ins_E$ in much the same way as the binary edge-insertion in Section 3.1 of \cite{xiao2026hopf}; the only difference is that trees are no longer binary.

Let $T, S \in \tilde \cT$ and let $e$ be an edge of $T$. To insert $S$ into $T$ at $e$, split the edge $e$ to create a new vertex and join the root of $S$ to this vertex by a new edge. We define $T \ins_E S$ as the sum of trees produced by inserting $S$ at each of edge of $T$, plus the tree obtained by grafting $T$ and $S$ to a common root, which coincides with $\fM(T,S)$.

\begin{example}
   Let $T = \tinytree{[[][][]]}$ and $S = \bullet$. Then
  \[
    T \ins_E S
    =
    3 \tinytree{[[[][]][][]]}
    +
    \tinytree{[[[][][]][]]}.
  \]
\end{example}

\begin{proposition}
  $\ins_E$ is a pre-Lie operation.
\end{proposition}

\begin{proof}
  The proof is identical to the proof of Theorem 3.1 from \cite{xiao2026hopf}. Indeed, nowhere in that proof did we use the fact that trees are binary. We repeat it here for convenience.


    Let $T_1, T_2, T_3$ be three trees. All trees are viewed as having a root edge. We will assign a geometric meaning to the associator,
  \[ A(T_1,T_2,T_3) = (T_1 \ins T_2) \ins T_3 - T_1 \ins (T_2 \ins T_3), \]
  which will show that the associator is unchanged when the positions of $T_2$ and $T_3$ are swapped.
  
  Note that $T_2 \ins T_3$ contains the term $\fM(T_2,T_3)$: the terms in $T_1 \ins (T_2 \ins T_3)$ resulting from $T_1 \ins \fM(T_2,T_3)$ correspond exactly to the terms in $(T_1 \ins T_2) \ins T_3$ where $T_2$ is inserted into an edge of $T_1$, and then $T_3$ is inserted into the edge which was originally the root of $T_2$.

  The terms in the associator $A(T_1,T_2,T_3)$ consist of the trees obtained when $T_2$ and $T_3$ either are inserted into distinct edges of $T_1$, or are both inserted into the same edge. In the latter case, each edge produces two terms: upon inserting $T_2$ into $T_1$ at an edge $e$, the edge is split, so $T_3$ is subsequently inserted into either the upper or lower edge, so (the root of) $T_3$ appears either above or below (the root of) $T_2$.

  The exact same terms appear in $A(T_1,T_2,T_3)$, since the order in which $T_2$ and $T_3$ are inserted is ultimately irrelevant. Hence $\ins$ satisfies the Vinberg identity, making it a pre-Lie operator.
\end{proof}

This insertion only introduces new internal vertices with binary branching. Let $\cT_{bin} \subset \tilde \cT$ be the collection of reduced binary trees, which includes the single vertex tree. When $\ins_E$ is restricted to the vector space $k\{ \cT_{bin} \}$, its image is closed in $k\{ \cT_{bin} \}$, meaning that $(k\{ \cT_{bin} \}, \ins_E)$ is a pre-Lie subalgebra of $(k\{ \tilde \cT \}, \ins_E)$. Then the Hopf algebra of $\cH_{bin}$ in \cite{xiao2026hopf} is recovered by taking the dual of $\cU(k\{ \cT_{bin} \}, \ins_E)$.

\begin{remark}
  When considering edge-insertion of reduced trees, it is often convenient to view each tree as having a virtual edge at the root. Insertion into the virtual edge corresponds to the Merge term. Then $T \ins_E S$ is a sum indexed by all the edges of $T$, including the virtual edge.
\end{remark}

\subsection{Combining vertex- and edge-insertion on reduced trees}

We can define another pre-Lie insertion operator $\ins_{red} = \ins_\ell + \ins_E$ on reduced trees which allows insertion at both internal vertices and edges. While this appears to be an unholy amalgamation of the operations described in the previous two sections, it is fairly natural that these two concepts can be viewed as part of the same operation when working with reduced trees.

\begin{example}
  We provide some small examples of $\ins_{red}$ below:
  \begin{align*}
    \bullet \ins_{red} \bullet &= \tinytree{[[][]]} \\
    \tinytree{[[][][]]} \ins_{red} \bullet
    &=
    \tinytree{[[][][][]]}
    +
    3 \tinytree{[[[][]][][]]}
    +
    \tinytree{[[[][][]][]]}
  \end{align*}
\end{example}

\begin{proposition}
  $\ins_{red}$ is a pre-Lie operation.
\end{proposition}

\begin{proof}
  Let $T_1, T_2, T_3$ be three reduced trees. As with the case of vertex-insertion and edge-insertion, the associator $A_{\ins_{red}}(T_1, T_2, T_3)$ captures ``independent'' insertions of $T_2$ and $T_3$ into $T_1$. The independent insertions fall into the following categories:
  \begin{enumerate}[(i)]
    \item $T_2$ and $T_3$ are inserted into the same vertex,
    \item $T_2$ and $T_3$ are into inserted different locations on the same edge,
    \item $T_2$ and $T_3$ are inserted into the same location on the same edge,
    \item or $T_2$ and $T_3$ are inserted into different locations (edge or vertex).
  \end{enumerate}
  Note that (i) is the case of vertex-insertion in the sense of $\ins_\ell$ and (ii) is the case of edge insertion in the sense of $\ins_E$.  In all cases, the same terms appear in $A_\ins(T_1, T_3, T_2)$, so the associators are equal and $\ins_{red}$ is a pre-Lie operation.
\end{proof}

\begin{remark}
  Lo\"{i}c Foissy points out that there is a stronger notion of ``compatibility'' between $\ins_\ell$ and $\ins_E$, as pre-Lie operators, on $k\{\tilde\cT\}$. Let $T_1, T_2, T_3$ be three reduced rooted trees. Then
  \begin{align*}
    (T_1 \ins_\ell T_2) &\ins_E T_3 - T_1 \ins_\ell (T_2 \ins_E T_3) - (T_1 \ins_E T_2) \ins_\ell T_3 + T_1 \ins_E (T_2 \ins_\ell T_3) \\
    &=
    (T_1 \ins_\ell T_3) \ins_E T_2 - T_1 \ins_\ell (T_3 \ins_E T_2) - (T_1 \ins_E T_3) \ins_\ell T_2 + T_1 \ins_E (T_3 \ins_\ell T_2)
  \end{align*}
  The operators $\ins_\ell$ and $\ins_E$ are compatible with each other in the sense of Strohmayer \cite{strohmayer2008operads}. In fact, $a \ins_E + b \ins_\ell$ is a pre-Lie operator for any scalars $a,b \in k$.
\end{remark}

Any reduced tree $T \in \tilde \cT$ admits a lift $\lambda(T) \in \cT$ where every edge of $T$ is subdivided and an additional edge is appended to the root. The subdivision of each edge creates an insertion spot at each edge and the new root motivates the presence of the Merge term. For example,
\[\lambda(\bullet) = \tinytree{[[]]}, \quad
 \lambda\left( \tinytree{[[[][]][][]]} \right)
 =
 \tinytree{[[[[[[]][[]]]][[]][[]]]]}.
\]
The terms of $T \ins_{red} S$ are then in correspondence with the terms of $\lambda(T) \ins_\ell \lambda(S)$. Indeed, if $\pi : k\{\cT\} \to k\{ \tilde \cT \}$ is the projection map such that $\pi \circ \lambda = \id$, then
\[ T \ins_{red} S = \pi(\lambda(T) \ins_\ell \lambda(S)). \]

Although this relationship holds at the level of trees, it does not automatically hold when extending to forests. To allow multiple insertions into the same edge, not necessarily at the same vertex, we would require the lift of that edge in $\cT_{red}$ to be broken up into as many vertices as there are trees to insert.

After an edge insertion splits an edge and introduces a new vertex, additional trees can be inserted into the new vertex in the extension to $k[\cT_{red}]$. Note that it does not matter which tree was inserted first, since trees are nonplanar and commutative. For intuition, it is helpful to examine the associator in a small example:

\begin{example}
  Let $T = \tinytree{[[][]]}$ and $S_1 = {\color{blue}\bullet}, S_2 = {\color{purple}\bullet}$. For illustrative purposes, we pretend that the colors are distinct and depict the new edges created by insertion with dashed lines. Then we will combine the terms which are the same isomorphism class of tree. The associator is
  \begin{align*}
    A_{\ins_{red}}(T, S_1, S_2) =
    \tinytree{[[][][,edge=dashed,color=blue][,edge=dashed,color=purple]]}
    &+ \tinytree{[[[][,edge=dashed,color=purple]][][,edge=dashed,color=blue]]}
    + \tinytree{[[[][,edge=dashed,color=blue]][][,edge=dashed,color=purple]]}
    + \tinytree{[[[][,edge=dashed,color=blue][,edge=dashed,color=purple]][]]}
    + \tinytree{[[[[][,edge=dashed,color=purple]][,edge=dashed,color=blue]][]]}
    + \tinytree{[[[[][,edge=dashed,color=blue]][,edge=dashed,color=purple]][]]} \\
    &+ \tinytree{[[][[][,edge=dashed,color=purple]][,edge=dashed,color=blue]]}
    + \tinytree{[[][[][,edge=dashed,color=blue]][,edge=dashed,color=purple]]}
    + \tinytree{[[][[][,edge=dashed,color=blue][,edge=dashed,color=purple]]]}
    + \tinytree{[[][[[][,edge=dashed,color=purple]][,edge=dashed,color=blue]]]}
    + \tinytree{[[][[[][,edge=dashed,color=blue]][,edge=dashed,color=purple]]]} \\
    = \tinytree{[[][][][]]} &+ 4 \tinytree{[[[][]][][]]} + 2 \tinytree{[[[][][]][]]} + 4 \tinytree{[[[][]][]]}.
  \end{align*}
\end{example}

\begin{remark}
$\ins_{red}$ is not free because
\[
  \tinytree{[[][][]]}
  =
  \tinytree{[[][]]} \ins \bullet
  -
  3 \left( \bullet \ins \tinytree{[[][]]} \right).
\]
\end{remark}

\section{Hopf algebras of trees}
\label{trees}

We remind the reader that $\cT$ is the collection of unlabelled rooted nonplanar trees, while $\tilde \cT$ is the collection of reduced trees.

\subsection{Symmetries of forests}
\label{symmetry}

Let $F$ be a forest of trees, either all in $\cT$ or all in $\tilde \cT$. Let $SG(F)$ be the group of automorphisms on $F$, which are graph isomorphisms that map roots to roots and that preserve the implicit directed graph structure induced by the root node.

We write $s_F = |SG(F)|$ for the cardinality of the automorphism group. Each graph automorphism induces an automorphism on the edges and vertices of $F$.

If $F = T_1^{a_1} \cdots T_k^{a_k}$, where $T_1, \dots, T_k$ are all distinct trees and each $a_i$ is a positive integer exponent, then
\[ s_F = a_1 ! \cdots a_k! \cdot s_{T_1}^{a_1} \cdots s_{T_k}^{a_k}, \]
which follows directly from the identification
\[ SG(F) \cong \fS_{a_1} \oplus \cdots \oplus \fS_{a_k} \oplus SG(T_1)^{\oplus a_1} \oplus \cdots \oplus SG(T_k)^{\oplus a_k}. \tag{$\star$} \]

\begin{example}
We provide the values of $s_T$ for some small trees below:

\begin{center}
\begin{tabular}{|c||c|c|c|c|c|c|c|c|c|c|c|c|c|}
\hline
$T$ & $\bullet$
    & $\tinytree{[[][]]}$
    & $\tinytree{[[[][]]]}$
    & $\tinytree{[[[]][]]}$
    & $\tinytree{[[[[]]][[]][]]}$
    & $\tinytree{[[][][]]}$
    & $\tinytree{[[[][]][]]}$
    & $\tinytree{[[][][][]]}$ 
    & $\tinytree{[[[][][]][]]}$ 
    & $\tinytree{[[[][]][][]]}$ 
    & $\tinytree{[[[][]][[][]]]}$ 
    & $\tinytree{[[[[][]][]][]]}$
    \\ \hline\hline
$s_T$ & 1 & 2 & 2 & 1 & 1 & 6 & 2 & 24 & 6 & 4 & 8 & 2 \\ \hline
\end{tabular}
\end{center}
\end{example}

We now define an inner product on $k[\cT]$ (or, equivalently, on $k[\tilde\cT]$) by
\[ \langle F, F' \rangle := s_F \delta_{F,F'} \]
for any monomials $F,F'$ in $k[\cT]$ (or $k[\tilde\cT]$). This pairing will play a key role in our definition of duality of Hopf algebras of forests.

\begin{definition}
  Bruned's Hopf algebra $(\cH_B, \sqcup, \Delta_B)$ is defined on $k[\cT]$, with a grading based on the number of leaves in a forest. The product is given by disjoint unions of forests. The coproduct is defined as
  \[ \Delta_B(F) = \sum_{a+b =n} (p_a \otimes p_b) \Delta_{CK}(F), \]
  where $\Delta_{CK}$ is the coproduct of $\cH_{CK}$ based on admissible cuts. That is, the coproduct $\Delta_B$ retains the terms of $\Delta(F)$ which do not introduce any additional leaves upon cutting.
\end{definition}

\begin{example}
Let $T = \tinytree{[[[[]]][[]][]]}$. Then
\[
  \Delta_B(T)
  =
  1 \otimes T + T \otimes 1
  + \tinytree{[[[]]]} \otimes \tinytree{[[[]][]]}
  + \tinytree{[[]]} \otimes \tinytree{[[[[]]][]]}
  + \bullet \otimes \tinytree{[[[[]]][[]]]}
  + \tinytree{[[[]]]} \tinytree{[[]]} \otimes \tinytree{[[]]}
  + \tinytree{[[[]]]} \bullet \otimes \tinytree{[[[]]]}
  + \tinytree{[[]]} \bullet \otimes \tinytree{[[[[]]]]}.
\]
\end{example}

\begin{remark}
  The terms in the coproduct $\Delta_B(F)$  are in one-to-one correspondence with subsets of leaves of $F$, just as in the binary case.
\end{remark}

That $\Delta'$ is coassociative will immediately follow from the next proposition.

\begin{proposition}
  The universal enveloping algebra $\cU(k\{\cT\}, \ins_\ell)$ is dual to the opposite of Bruned's Hopf algebra $\cH_B$.
\end{proposition}

\begin{proof}
  The duality between $\cU(k\{\cT\}, \ins_{free})$ and (the opposite of) $\cH_{CK}$ with respect to $\langle -, -\rangle$ is known by \cite{hoffman2003combinatorics}, and $\ins_\ell$ is a restriction of $\ins_{free}$ induced by the leaf-based grading, so $\Delta_B$ is exactly the coproduct described in part (3) of Theorem \ref{compat}.
\end{proof}

Let $\pi : k[\cT] \to k[\tilde \cT]$ be the map which takes a tree $T$ to its reduced equivalent; that is, all non-root vertices of degree $2$ are flattened into their neighboring edges. For example,

\[ \pi \left( \tinytree{[[[[]]][[]][]]} \right) = \tinytree{[[][][]]}. \]

\begin{definition}
  The Hopf algebra of reduced forests $(\cH_{red}, \sqcup, \Delta_{red})$ is defined on $k[\tilde\cT]$ with a grading based on the number of leaves in a forest. The product is given by disjoint unions of forests. The coproduct is defined as
  \[ \Delta_{red}(\pi(F)) = (\pi \otimes \pi)(\Delta_B(F)) \]
  where $F$ is any lift of $\pi(F) = \tilde F \in k[\cT]$.
\end{definition}

\begin{remark}
  It immediately follows from this definition that $\pi : \cH_B \to \cH_{red}$ is a Hopf algebra morphism, and $\cH_{red}$ is a quotient of $\cH_B$.
\end{remark}

Foissy studies both $\cH_B$ and $\cH_{red}$ extensively in \cite{foissy2026combinatorial} from the perspective of combinatorial twisted bialgebras, showing also that the corresponding coproducts are well-defined. Our alternative proof of the coassociativity of $\Delta_B$ highlights how $(\cH_B, \Delta_B)$ can be derived by restricting a pre-Lie operation, taking the universal enveloping algebra, and then invoking duality. Foissy also provides cancellation-free formulas for the antipode of both Hopf algebras using the quasi-commutativity of these spaces as combinatorial twisted bialgebras.

\begin{example}
Let $T = \tinytree{[[][][]]}$. Then
\[
  \Delta_{red}(T)
  =
  1 \otimes T + T \otimes 1
  + 3 \bullet \otimes \tinytree{[[][]]}
  + 2 \bullet\bullet \otimes \bullet.
\]
\end{example}

That $\cH_{red}$ is related to $\ins_{red}$, the insertion operation on reduced trees involving insertions at both edges and internal vertices, requires a proof. This essentially combines arguments from the duality of vertex-insertion from \cite{hoffman2003combinatorics} and duality of edge-insertion from \cite{xiao2026hopf}.

\begin{remark}
  $\cH_{bin}$ is a Hopf subalgebra of $\cH_{red}$.
\end{remark}

\begin{theorem}
  The Hopf algebra $\cH_{red}$ is dual to the opposite of $\cU(k\{\tilde \cT\}, \ins_{red})$, in the sense that
  \[ \langle F_1 \ast F_2, F_3 \rangle = \langle F_2 \otimes F_1, F_3 \rangle.  \]
\end{theorem}

\begin{proof}
  We write $\cU(\cL) = \cU(k\{\tilde \cT\}, \ins_{red})$ to simplify notation.

  The duality between the product $\sqcup$  of $\cH_{red}$ and the coproduct $\Sha$ of $\cU(\cL)$ is clear. Indeed, the proof is the same for any polynomial Hopf algebra and corresponding unshuffle coproduct; see Section 4.3 of \cite{xiao2026hopf}. We will follow the notation used in both \cite{hoffman2003combinatorics} and \cite{xiao2026hopf}, defining multiplicities of cuts and grafts for reduced forests. The ultimate goal is to show that
  \[ \langle F_1 \ast F_2, F_3 \rangle = \langle F_2 \otimes F_1, F_3 \rangle.  \]

  It is convenient to view all reduced trees as having a virtual root edge. When speaking of the edge set $E(F)$, the virtual root edge of each tree in $F$ is included.

  Given a forest $F \in k[\tilde \cT]$, let $\cC(F)$ be the collection of leaf-admissible cuts on $F$. Each $C$ in $\cC$ is a subset of $E(F)$ such that
  \begin{enumerate}
    \item $C$ is an admissible cut in the sense of Connes-Kreimer, in that no pair of edges in $C$ appears on the same path from root to leaf of any tree;
    \item for no vertex do all of its child-edges appear in $C$.
  \end{enumerate}
  The second criterion explicitly excludes the creation of additional leaves upon cutting.

  We now write $\Delta = \Delta_{red}$ to simplify notation. The coproduct is now a sum over cuts in~$\cC(F)$:
  \[ \Delta(F) = \sum_{C \in \cC(F)} P^C(F) \otimes R^C(F), \]
  where $P^C(F)$ is the forest consisting of the trees lying below the cut edges, and $R^C(F)$ is the reduced remainder. If a root edge appears in $C$, then the tree below it appears in $P^C(F)$ and has no remainder in $R^C(F)$.

  To enumerate grafts, consider two reduced forests $F_1$ and $F_2$ and let $\cG(F_1, F_2)$ be the collection of distinct grafts of $F_1$ into $F_2$. When speaking of grafts, we include the terms introduced by the extension of $\ins$ to $\ast$. Each distinct graft contributes a forest $F_2 \ast_G F_1$ to $F_2 \ast F_1$, so
  \[ F_2 \ast F_1 = \sum_{G \in \cG(F_1, F_2)} F_2 \ast_G F_1. \]
  Each graft consists of the following data:
  \begin{enumerate}
    \item a map $\mu_G$ from the trees in $F_1$ to $E(F_2) \sqcup IV(F_2) \sqcup \{e_0\}$, where $IV(F_2) \subseteq V(F_2)$ are the internal vertices of $F_2$, and $e_0$ is the image of any tree not inserted into an edge or vertex of $F_2$;
    \item a composition (ordered partition) on $\mu_G^{-1}(e)$ for every $e \in E(F_2)$. When two trees are in the same component of the partition, that means any one of the trees is first grafted into the edge, creating a new vertex, and then additional trees are grafted into the vertex that has been created. It is not relevant which tree is chosen as the ``first'' graft, since this insertion is simultaneous.
    
  \end{enumerate}
  No composition or any linear ordering is imposed on the preimage of $e_0$ or any vertex in $IV(F_2)$, since trees commute and the children of any vertex are unordered due to non-planarity of trees.

  Recall that in the case of edge-insertion in \cite{xiao2026hopf}, a linear order was imposed on $\mu_G^{-1}(e)$ since each instance of edge-insertion creates a new vertex. A linear order is essentially a composition on $\mu_G^{-1}(e)$ where each component is a singleton, so this is a natural generalization.

  Now we have established the notions of cuts and grafts for reduced forests. Let $F_1, F_2$ and $F_3$ be reduced forests. We define
\begin{align*}
  n(F_1, F_2; F_3)  &= | \{ G \in \cG(F_1, F_2) : F_2 \ast_G F_1 = F_3 \}|  \\
  &= \text{the coefficient of $F_3$ in $F_2 \ast F_1$}
\end{align*}
and
\begin{align*}
  m(F_1, F_2; F_3) &= | \{ C \in \cC(F_3) : P^C(F_3) = F_1 \text{ and } R^C(F_3) = F_2 \} | \\
  &= \text{the coefficient of $F_1 \otimes F_2$ in $\Delta(F_3)$},
\end{align*}
so that
\[ \langle F_1 \ast F_2, F_3 \rangle = n(F_1, F_2; F_3) \langle F_3, F_3 \rangle = n(F_1, F_2; F_3) s_{F_3} \]
and
\begin{align*}
  \langle F_2 \otimes F_1, \Delta(F_3) \rangle
  &= m(F_1, F_2; F_3) \langle F_2 \otimes F_1, F_2 \otimes F_1 \rangle \\
  &= m(F_1, F_2; F_3) \langle F_2, F_2 \rangle \langle F_1, F_1 \rangle \\
   &= m(F_1, F_2; F_3) s_{F_1} s_{F_2}. \end{align*}
 Having fixed $F_1, F_2, F_3$, our aim is to show that
 \[ n(F_1, F_2; F_3) s_{F_3} = m(F_1, F_2; F_3) s_{F_1} s_{F_2}. \]
 Let
 \[ \cC' = \cC_{F_1, F_2, F_3}' := \{ C \in \cC(F_3) : P^C(F_3) = F_1, R^C(F_3) = F_2 \}, \]
 so that $|\cC'| = n(F_1, F_2; F_3)$.
 
As with the case of binary trees and edge insertion in \cite{xiao2026hopf}, $\cC'$ has a decomposition as a disjoint union of \emph{cut-classes}, each of which is the orbit of any cut in its class under the action of $SG(F_3)$ on $E(F_3)$.

With the same forests $F_1, F_2, F_3$ as above, we define
\[ \cG' = \{ G \in \cG(F_1, F_2) : F_2 \ast_G F_1 = F_3 \} \]
so that $m(F_1, F_2; F_3) = |\cG'|$, and we partition $\cG'$ into graft-classes $\cG''$ under the action of $SG(F_2)$ on $E(F_2) \sqcup IV(F_2)$, allowing grafting at internal vertices as well as edges.

We allow the reader to appeal to \cite{xiao2026hopf} for the details of the cut-classes and graft-classes and the application of the orbit-stabilizer theorem for the case of edge-insertion, and to convince themself that the proof remains fundamentally unchanged when allowing insertion at internal vertices.

The only difference is that when calculating the cardinality of $\cG'$ as a sum over the cardinalities of graft classes $|\cG''|$, certain vertex insertions are simultaneous, so instead of each $|\cG''|$ being multiplied by the factorials of the number of copies of each tree in $F_3$, which is the same for every graft-class, multiple copies of the same tree may belong to the same simultaneous insertion. In that case they are interchangeable, so we divide the coefficient by the number of identical trees in the same set in each graft-class.
\end{proof}

The following example should aid the reader's intuition when considering cut-classes and graft-classes.

\begin{example}
We illustrate an example where
\[ F_1 = \tinytree{[[][]]}, \quad F_2 = \bullet \bullet\bullet \tinytree{[[][]]}, \quad F_3 = \tinytree{[[[[[][]][]][[][]]][][]]}, \] 
which have symmetry coefficients
\[ s_{F_1} = 2, s_{F_2} = 12, s_{F_3} = 8.\]
In the middle column of the following table, whenever we label an edge with a composition, the parts of the composition when read left-to-right are inserted top-to-bottom.

We include the virtual root edge in the illustrations of the grafts, but not in the illustrations of the cuts.

\forestset{
    bcut/.style={edge path={\noexpand\path[\forestoption{edge}](!u.parent anchor)-- node[red] {\scalebox{1}{$\diagdown$}} (.child anchor)\forestoption{edge label};}},
    fcut/.style={edge path={\noexpand\path[\forestoption{edge}](!u.parent anchor)-- node[red] {\scalebox{1}{$\diagup$}} (.child anchor)\forestoption{edge label};}},
    scut/.style={edge path={\noexpand\path[\forestoption{edge}](!u.parent anchor)-- node[red] {\scalebox{1.5}{$-$}} (.child anchor)\forestoption{edge label};}},
}

\begin{center}
\begin{tabular}{ |c||c||c| }
     \hline
Cut-classes in $\cC'$ & Graft-classes in $\cG'$ & Coeff. of $|\cG''|$ \\
     \hline\hline

\(
\bigtree{[[[[,bcut[][]][]][[,bcut][]]][,scut][,fcut]]},
\bigtree{[[[[,bcut[][]][]][[][,fcut]]][,scut][,fcut]]}
\)
&

\(
\biggertree{[
[,edge={blue,very thick},edge label={node[blue,midway,right,font=\scriptsize]{ $(\{\bullet,\bullet\})$ }}
[,edge={blue,very thick},edge label={node[blue,midway,left,font=\scriptsize]{ $(\{\wedge)\}$ }}]
[,edge={blue,very thick},edge label={node[blue,midway,right,font=\scriptsize]{ $(\{\bullet\})$ }}
]
]]},
\biggertree{[
[,edge={blue,very thick},edge label={node[blue,midway,right,font=\scriptsize]{ $(\{\bullet,\bullet\})$ }}
[,edge={blue,very thick},edge label={node[blue,midway,left,font=\scriptsize]{ $(\{\bullet\})$ }}]
[,edge={blue,very thick},edge label={node[blue,midway,right,font=\scriptsize]{ $(\{\wedge\})$ }}
]
]]}
\)
&
3

\\ \hline

\(
\bigtree{[[[[[][]][,fcut]][,fcut[][]]][,scut][,fcut]]}
\)
&

\(
\biggertree{[
[,edge={blue,very thick},edge label={node[blue,midway,right,font=\scriptsize]{ $(\{\bullet,\bullet\}, \{\wedge\}, \{\bullet\})$ }}
[]
[]
]]}
\)

&
3

\\ \hline

\(
\bigtree{[[[[[,bcut][]][]][,fcut[][]]][,scut][,fcut]]},
\bigtree{[[[[[][,fcut]][]][,fcut[][]]][,scut][,fcut]]}
\)
&

\(
\biggertree{[
[,edge={blue,very thick},edge label={node[blue,midway,right,font=\scriptsize]{ $(\{\bullet,\bullet\}, \{\wedge\})$ }}
[,edge={blue,very thick},edge label={node[blue,midway,left,font=\scriptsize]{ $(\{\bullet\})$ }}]
[]
]]},
\biggertree{[
[,edge={blue,very thick},edge label={node[blue,midway,right,font=\scriptsize]{ $(\{\bullet,\bullet\}, \{\wedge\})$ }}
[]
[,edge={blue,very thick},edge label={node[blue,midway,right,font=\scriptsize]{ $(\{\bullet\})$ }}]
]]}
\)
&
3

\\ \hline

\(
\begin{aligned}
\\[-3mm]
\bigtree{[[[[[,bcut][]][,fcut]][,fcut[][]]][,scut][]]},
\bigtree{[[[[[][,fcut]][,fcut]][,fcut[][]]][,scut][]]},
\\
\bigtree{[[[[[,bcut][]][,fcut]][,fcut[][]]][][,fcut]]},
\bigtree{[[[[[][,fcut]][,fcut]][,fcut[][]]][][,fcut]]}
\\[-3mm]
\end{aligned}
\)

&

\(
\biggertree{[
[,color={purple},
edge label={node[purple,left,font=\scriptsize]{ $\{\bullet\}$ }}
[
,edge={blue,very thick},
edge label={node[blue,midway,left,font=\scriptsize]{ $(\{\wedge\}, \{\bullet\}, \{\bullet\})$ }}]
[]
]]},
\biggertree{[
[,color={purple},
edge label={node[purple,right,font=\scriptsize]{ $\{\bullet\}$ }}
[]
[,edge={blue,very thick},
edge label={node[blue,midway,right,font=\scriptsize]{ $(\{\wedge\}, \{\bullet\}, \{\bullet\})$ }}]
]]}
\)
&
6

\\ \hline
\end{tabular}
\end{center}

We observe that the coefficients in the first three rows are $3 = \frac{3!}{2!}$, due to the appearance of the multiset $\{\bullet, \bullet\}$. On the other hand, each copy of $\bullet$ is in its own set in the last row, so the coefficient is $6 = 3!$. We count all the elements in each cut-class, which gives
\[ |\cC'| =  2 + 1 + 2 + 4 = 9. \]
The weighted sum of graft-classes is
\[ |\cG'| = 3\cdot 2 + 3 \cdot 1 + 3 \cdot 2 + 6 \cdot 2 = 27. \]
Finally,
\[ |\cC'| s_{F_1} s_{F_2} = 9 \cdot 2 \cdot 12 = 216 = 27 \cdot 8 = |\cG'| s_{F_3}.  \]

\end{example}

\subsection{The pre-Lie exponential}

In Section 5 of \cite{xiao2026hopf}, we derived some consequences of the dual pairing for binary trees. The results can be easily extended to $\cH_{red}$, since the proofs rely only on the same duality extended to reduced trees which are not necessarily binary. We specialize the multiplicities $n(F_1, F_2; F_3)$ and $m(F_1, F_2; F_3)$ to the case where $F_2 = T$ is a tree and $F_1 = \bullet$. Let $\tilde \cT_n$ be the collection of reduced trees with $n$ leaves.

Following Hoffman's notation, define the growth operator $\fN : k\{\cT_n\} \to k \{\cT_{n+1}\}$ by
\[ \fN(T) = T \ins \bullet = \sum_{T' \in \cT} n(\bullet, T; T') T' \]
and the pruning operator $\fP : k\{\cT_n\} \to k\{\cT_{n-1}\}$ by
\[ \fP(T) = \sum_{T' \in \cT} m(\bullet, T'; T) T', \]
which are adjoint, by Proposition 5.1 of \cite{xiao2026hopf}.

The pre-Lie exponential of an element $x \in \cL$ is defined as the formal series
\[ W(x) = x + \frac{1}{2!} x \ins x + \frac{1}{3!} (x \ins x) \ins x + \cdots. \]
When equipped with a multiplication operation
\[ W(x) \star W(y) = W(C(x,y)) \]
in terms of the Baker-Campbell-Hausdorff formula
\[ C(x,y) = x + y + \frac{1}{2} [x,y] + \frac{1}{12} ( [x,[x,y]] + [y,[y,x]] ) - \frac{1}{24} [y,[x,[x,y]]] + \cdots, \]
the formal series $W(x)$ have the structure of a multiplicative group. See \cite{marcolli2015graph} for more references and Section 3.6 of \cite{cartier2021classical} for the formula of the full expansion of the Baker-Campbell-Hausdorff formula.

Let $\bullet \in \cH_{red}$. Then
\[ W(\bullet) = \sum_{k=1}^\infty \frac{1}{k!} \fN^{k-1}(\bullet), \]
where $\fN^0(\bullet) = \bullet$.

Just as with binary trees and edge-insertion, the coefficient of each tree in $W(\bullet)$ has a combinatorial interpretation, which follows immediately from duality.

The first few iterated insertions of $\bullet$ are:

  \begin{align*}
    \fN(\bullet) &= \tinytree{[[][]]} \\
    \fN^2(\bullet) &= \tinytree{[[][][]]} + 3 \tinytree{[[[][]][]]} \\
    \fN^3(\bullet) &=
    \tinytree{[[][][][]]}
     + 4 \tinytree{[[[][][]][]]}
     + 6 \tinytree{[[[][]][][]]}
     + 3 \tinytree{[[[][]][[][]]]}
     + 12 \tinytree{[[[[][]][]][]]}
  \end{align*}

Thus, the terms in $W(\bullet)$ with up to four leaves are
  \begin{align*}
    W(\bullet) =
    \bullet &+ \frac{1}{2} \tinytree{[[][]]} 
    + \frac{1}{6} \tinytree{[[][][]]} + \frac{1}{2} \tinytree{[[[][]][]]} \\
    &+ \frac{1}{24} \tinytree{[[][][][]]}
     + \frac{1}{6} \tinytree{[[[][][]][]]}
     + \frac{1}{4} \tinytree{[[[][]][][]]}
     + \frac{1}{8} \tinytree{[[[][]][[][]]]}
     + \frac{1}{2} \tinytree{[[[[][]][]][]]} + \cdots
  \end{align*}

\begin{proposition}
  The coefficient of $T$ in $W(\bullet)$ is $1/s_T$.
\end{proposition}

\begin{proof}
  The proof is identical to that of Proposition 5.2 of \cite{xiao2026hopf}.
\end{proof}

\section{Hochschild cohomology}
\label{cocycle}

We refer to Chapter 2.4 of Erik Panzer's MSc thesis \cite{panzer2012msc} as our primary reference for Hochschild cohomology.

\begin{definition}
  Let $\cM$ be a vector space and $(\cC, \Delta, \eps)$ a coalgebra. Then $\cM$ has a \emph{right $\cC$-comodule} structure if there is a linear map $\rho : \cM \to \cM \otimes \cC$ satisfying
  \[ (\id_\cM \otimes \Delta)\circ \rho = (\rho \otimes \id_\cC) \circ \rho \]
  and
  \[ \id_\cM = (\id_\cM \otimes \eps) \circ \rho. \]

  Analogously, $\cM$ is a \emph{left $\cC$-comodule} if there is a linear map
  $\lambda : \cM \to \cM \otimes \cC$
  such that
  \[ (\Delta \circ \id_\cM) \circ \lambda = (\id_\cC \otimes \lambda) \circ \lambda \]
  and
  \[ \id_\cM = (\eps \otimes \id_\cM) \circ \lambda. \]

  If $\cM$ has both a left and right $\cC$-comodule structure, then we call $\cM$ a \emph{$\cC$-bicomodule}.
\end{definition}

\begin{definition}
  Let $(\cM, \Delta)$ be a vector space. Suppose $\cM$ has a $\cC$-bicomodule structure with left coaction $\lambda : \cM \to \cC \otimes \cM$ and right coaction $\rho : \cM \to \cM \otimes \cC$. The \emph{Hochschild cochain complex} is defined as 
  \[ C^n(\cC, \cM) = \Hom(\cM, \cC^{\otimes n}) \]
  with differential
  \[ d^n\phi(x) = (\id \otimes \phi)\lambda(x) + \sum_{i=1}^n (-1)^i (\id^{\otimes (i-1)} \otimes \Delta \otimes \id^{\otimes(n-i)})\phi(x) + (-1)^{n+1} (\phi \otimes \id) \rho(x) \] 
  for $\phi  \in C^n(\cC,\cM)$. The $n$th Hochschild cohomology group is defined by
  \[ HH^n(\cC, \cM) := \frac{\ker d^n}{\im d^{n-1}}. \]
\end{definition}

We are primarily interested in the case where $\cM = \cC = \cH$, where $\cH$ is a Hopf algebra that is viewed as bicomodule over itself with left coaction $\lambda(h) = \Delta(h)$ and right coaction $\rho(h) = h \otimes 1$.

We will abbreviate $HH^\bullet(\cH,\cH)$ as $HH^\bullet(\cH)$, or even $H^\bullet(\cH)$ when there is no possibility for confusion (i.e. no other cohomology theories are considered).

Of particular interest are $1$-cocycles. A $1$-cochain $L : \cH \to \cH$ is a $1$-cocycle when it satisfies $d^1 L = 0$, that is
\[ d^1 L(h) = (\id \otimes L)\Delta(h) - \Delta L(h) + L(h) \otimes 1 = 0, \]
or equivalently,
\[ \Delta L(h) = (\id \otimes L)\Delta 
(h)  + L(h) \otimes 1 \tag{$\ast$} \]
for all $h \in \cH$. We will refer to ($\ast$) as the \emph{Hochschild $1$-cocycle condition}.

A map $A : \cH \to \cH$ is a $1$-coboundary when $A = d^0 \alpha$ for some $\alpha : \cH \to k$. In that case,
\[ A(h) = (d^0\alpha)(h) = (\id \otimes \alpha)\Delta(h) - \alpha(h) \otimes 1  \]
where the image of $A$ lies in $\cH$ from the canonical identification $\cH \cong \cH \otimes k \cong k \otimes \cH$. We note that in any unital Hopf algebra,
\[ (d^0\alpha)(1) = 1 \otimes \alpha(1) - \alpha(1) \otimes1 = 0 \]
so a sufficient condition for checking that a $1$-cocycle $L$ is nontrivial is to verify that $L(1) \neq 0$.

\subsection{Combining forests into a single tree}

We now define some forest-merging operations, which will be our main examples of Hochschild $1$-cocycles for the Hopf algebras of forests that we are interested in.

Define a linear map $B^+ : k[\cT] \to k[\cT]$ by
 \[
   B^+(1) = \bullet, \quad B^+(T_1 T_2 \cdots T_n) =
   {\vcenter{\hbox{
   \begin{forest}[,circle,fill,inner sep=0.6mm[$T_1$][$T_2$][$\cdots$][$T_n$]]\end{forest}
   }}},
 \]
where forests are joined into a single tree with a common root. The image is well-defined because monomials commute and trees are nonplanar. In particular, for a tree $T \in \cT$,
\[ B^+(T) = {\vcenter{\hbox{
   \begin{forest}[,circle,fill,inner sep=0.6mm[$T_1$]]\end{forest}
   }}} \]
so $B^+$ always increases the number of vertices by $1$.

A similar map $B^\circ : k[\cT] \to k[\cT]$ can be defined where $B^\circ(1) = 0$ and $B^\circ(F) = B^+(F)$ for non-empty forests.

Finally $B^\circ_{red}$ can be defined for forests of reduced trees. Here,
\[ B^\circ_{red}(1) = 0, \quad B^\circ_{red}(T) = T, \quad  B^\circ_{red}(T_1 T_2 \cdots T_n) =
   {\vcenter{\hbox{
   \begin{forest}[,circle,fill,inner sep=0.6mm[$T_1$][$T_2$][$\cdots$][$T_n$]]\end{forest}
   }}}. \]
The fact that $B^+$, $B^\circ$ and $B^\circ_{red}$ are examples of Hochschild 1-cocycles will justify the choices for the image of $1$ in each case.

\begin{remark}
  There is a simple description of the Connes-Kreimer $1$-cocycle $B^+$ in terms of the extension of $\ins_{free}$ to $k[\cT]$. For a forest $T_1 \sqcup \dots \sqcup T_n$,
  \[ B^+(T_1 \sqcup \cdots \sqcup T_n) = \bullet \ins_{free} (T_1 \sqcup \cdots \sqcup T_n), \]
  which is the simultaneous insertion of the forest into a single vertex.
\end{remark}

\begin{remark}
  In \cite{marcolli2026hochschild} a linear map $\K$ is defined on $\cH_{bin}$, the Hopf algebra of binary rooted forests from \cite{xiao2026hopf}. Following \cite{marcolli2026hochschild}, there is an alternative grading on $\cH_{bin}$ based on the number of connected components in a forest, so that $\cH^{(n)}_{bin} \subset \cH_{bin}$ is spanned by forests $F = T_1 \sqcup \cdots \sqcup T_n$ with $n$ components. We write $\Pi^{(n)} : \cH_{bin} \to \cH_{bin}^{(n)}$ for the projection map based on this grading.
  
  Then $\K : \cH_{bin} \to \cH_{bin}$ is defined as
  \[ \K := \Pi^{(1)} + B^+ \circ \Pi^{(2)}. \]
  where $B^+$ refers to the Connes-Kreimer $1$-cocycle. Although $\K$ is a not a Hochschild $1$-cocycle for $\cH_{bin}$, it satisfies the Hochschild $1$-cocycle condition for the $\cH_{bin}$-bicomodule
\[ \cC^{\leq 2} = \cH_{bin}^{(0)} \oplus \cH_{bin}^{(1)} \oplus \cH_{bin}^{(2)}, \]
which is the restriction of $\cH_{bin}$ to forests with at most $2$ components.

There is another definition of $\K$ based on the $1$-cocycle $B^\circ_{red}$. Let $ \iota : \cH_{bin} \hookrightarrow \cH_{red}$ be the inclusion map and let $\delta : \cH_{red} \twoheadrightarrow \cH_{bin}$ be the linear map which sends all forests which are not binary to $0$. Then
\[ \K = \delta \circ B^\circ_{red} \circ \iota. \]
\end{remark}

\begin{remark}
The permutation group, as defined in \ref{symmetry}, is invariant under $B^+$ for any non-empty forest,
\[ SG(F) \cong SG(B^+(F)), \]
so the permutation group of any tree can be computed inductively by examining the children of the root node and using ($\star$). The same holds for $B^\circ$ since the images of $B^\circ$ and $B^+$ coincide for forests.

If $F$ is a forest of reduced trees, then its permutation group is also invariant under $B^\circ_{red}$:
\[ SG(F) \cong SG(B^\circ_{red}(F)). \]
\end{remark}

\begin{example}
  $B^+$ is a nontrivial $1$-cocycle of $\cH_{CK}$; see \cite{connes1999hopf}, \cite{panzer2012msc}, \cite{foissy2013introduction} for proofs and details, including the universal property of the pair $(\cH_{CK}, B^+)$.
\end{example}

\begin{proposition}
  $B^\circ$ is a nontrivial $1$-cocycle of $\cH_B$.
\end{proposition}

\begin{proof}
  First we verify that $B^\circ$ satisfies the $1$-cocycle condition.
  
  The coproduct $\Delta_B$ of $\cH_B$ is the restriction of the coproduct $\Delta_{CK}$ on $\cH_{CK}$ which deletes all terms that create additional leaves. We recall the identity for a forest $F$ with $n$ leaves, with $p_k : \cH \to \cH_k$ being the projection to the degree $k$ component in the leaf-based grading:
  \[ \Delta_B(F) = \sum_{a+b=n} (p_a \otimes p_b) \Delta_{CK}(F). \]
  Since $B^\circ(F) = B^+(F)$ for non-empty forests (and both operations preserve leaf count for non-empty forests),
  \begin{align*}
    \Delta_B(B^\circ(F))
    &= \sum_{a+b=n} (p_a \otimes p_b) \Delta_{CK}(B^+(F)) \\
    &= \sum_{a+b=n} (p_a \otimes p_b) ((\id \otimes B^+)\Delta_{CK}(F) + B^+(F)\otimes 1) \\
    &= \sum_{a+b=n} (p_a \otimes p_b) (\id \otimes B^+)\Delta_{CK}(F) + B^\circ(F)\otimes 1
  \end{align*}
  where we have used the cocycle property for $B^+$ in the second equality. Let us examine which terms are retained in the sum after projection.
  
  Apart from $F \otimes B^+(1)$, every term $F' \otimes B^+(F'')$ in $(\id \otimes B^+)\Delta_{CK}(F)$ containing extraneous leaves in $B^+(F'')$  corresponds to a term $F' \otimes F''$ in $\Delta_{CK}(F)$ with extraneous leaves in $F''$. All such terms are eliminated by the projection.
  
  The term $F \otimes B^+(1)$ is also eliminated because $B^+(1) = \bullet$ has degree $1$, so we can identify $B^+$ with $B^\circ$ before projection.
  
  Furthermore, $p_a \otimes p_b$ and $\id \otimes B^\circ$ commute, so we can swap the positions of the two maps. That is,
  \begin{align*}
    \Delta_B(B^\circ(F))
    &= \sum_{a+b=n} (p_a \otimes p_b) (\id \otimes B^\circ)\Delta_{CK}(F) + B^\circ(F)\otimes 1 \\
    &= (\id \otimes B^\circ) \sum_{a+b=n} (p_a \otimes p_b) \Delta_{CK}(F) + B^\circ(F)\otimes 1 \\
    &= (\id \otimes B^\circ) \Delta_B(F) + B^\circ(F)\otimes 1.
  \end{align*}
  This shows that $B^\circ$ is a $1$-cocycle for $(\cH_B, \Delta_B)$.

  Next we check that $B^\circ$ is not a coboundary. Since $B^\circ(1) = 0$ it is not enough to verify that the image of $1$ is nonzero. Suppose, to the contrary, that $B^\circ = d\alpha$ for some $\alpha : \cH \to k$. Recalling, then, that
  \[ (d\alpha)(h) = (\id \otimes \alpha) \Delta_B(h) - \alpha(h) \otimes 1, \]
  we attempt to reconstruct such a map $\alpha$. We have
  \begin{align*}
    (d\alpha)(\bullet\bullet)
    &= (\id \otimes \alpha)(1 \otimes \bullet\bullet+2 \bullet\otimes\bullet+\bullet\bullet\otimes1) - \alpha(\bullet\bullet) \otimes 1 \\
    &= 1 \otimes \alpha(\bullet\bullet)
       + 2\bullet\otimes\alpha(\bullet)
       + \bullet\bullet \otimes \alpha(1)
       - \alpha(\bullet\bullet)\otimes 1 \\
       &= 2\alpha(\bullet) \bullet + \alpha(1) \bullet\bullet \\
       &= a \bullet + b \bullet\bullet
  \end{align*}
  for some coefficients $a,b \in k$. But $B^\circ(\bullet\bullet) = \tinytree{[[][]]}$, yet this term cannot appear in $(d\alpha)(\bullet\bullet)$. So $B^\circ$ is not a coboundary.
\end{proof}

\begin{remark}
  Since $B^\circ(F) = 0$ if and only if $F=1$, the terms of $\Delta_B B^\circ(F)$ are in bijective correspondence with the terms of $\Delta_B(F)$. Indeed, this is why Bruned was able to join all of the severed forests into a common root in his original presentation of $\cH_B$ in \cite{bruned2015singular}.
\end{remark}

\begin{proposition}
  $B^\circ_{red}$ is a nontrivial $1$-cocycle of $\cH_{red}$.
\end{proposition}

\begin{proof}
  Let $\bar F$ be a reduced forest in $k[ \tilde \cT ]$ and let $F$ be its lift to $k[\cT]$, which is also reduced. As a reminder, $\pi : k[\cT] \to k[\cT]$ is the projection map that sends any forest to a reduced forest, and it can also be viewed as a morphism of Hopf algebras $\pi : \cH_B \to \cH_{red}$.
  
  We aim to show that
  \[ \Delta_{red} B^\circ_{red}(\bar F) = (\id \otimes B^\circ_{red}) \Delta_{red} (\bar F) + B^\circ_{red}(\bar F) \otimes 1. \]
  Since $\pi$ is surjective, we can identify $\bar F$ with $\pi(F)$ and verify the equality for all $\pi(F)$. Also note that
  \[ B_{red}^\circ \circ \pi = \pi \circ B^\circ, \quad (\pi \otimes \pi) \circ \Delta_B = \Delta_{red} \circ \pi. \]
  The first identity shows that $B^\circ_{red}$ and $B^\circ$ are intertwined by the map $\pi$. (We will later show that $\pi$ satisfies a universal property related to the pairs $(\cH_B, B^\circ)$ and $(\cH_{red}, B_{red})$.) Then, using the fact that $B^\circ$ is a cocycle,
  \begin{align*}
    \Delta_{red} (B^\circ_{red}(\pi(F)))
    &= \Delta_{red} (\pi(B^\circ(F))) \\
    &= (\pi \otimes \pi) \Delta' (B^\circ(F)) \\
    &= (\pi \otimes \pi) ((\id \otimes B^\circ)\Delta'(F) + B^\circ(F) \otimes 1 ) \\
    &= (\pi \otimes \pi \circ B^\circ)\Delta'(F) + \pi(B^\circ(F) )\otimes 1 \\
    &= (\pi \otimes B^\circ_{red} \circ \pi)\Delta'(F) + B^\circ_{red}(\pi(F) )\otimes 1 \\
    &= (\id \otimes B^\circ_{red}) (\pi \otimes \pi)\Delta'(F) + B^\circ_{red}(\pi(F) )\otimes 1 \\
    &= (\id \otimes B^\circ_{red})\Delta_{red}(\pi(F)) + B^\circ_{red}(\pi(F) )\otimes 1.
  \end{align*}
  This shows that $B^\circ_{red}$ is a $1$-cocycle for $\cH_{red}$.
  
  To show that $B^\circ_{red}$ is not a coboundary, we can again appeal to the image of $\bullet \bullet$; the calculation is identical to the one in $\cH_B$ in the previous proof.
\end{proof}

Alternatively, we could have proved that $B^\circ_{red}$ is a $1$-cocycle using a combinatorial argument, bypassing the usage of the projection map $\pi$ to reduced forests. In fact, $\pi$ itself can be obtained from a universal property related to the pair $(\cH_B, B^\circ)$, which is similar to the universal property of $(\cH_{CK}, B^+)$ with the additional choice of an element in the target algebra. The fact that $\pi$ is a morphism of Hopf algebras also follows directly from the universal property, which is stated below.

\begin{theorem}
  If $\cA$ is a commutative unital algebra and $L : \cA \to \cA$ is a linear map such that $L(1) = 0$, and $a \in \cA$, then there is a unique algebra morphism $\phi_{L,a} : \cH_B \to \cA$ such that $\phi_{L,a}(\bullet) = a$ and
  \[ \phi_{L,a} \circ B^\circ = L \circ \phi_{L,a}. \]
  If $\cA$ is also a bialgebra, $a \in \Prim(\cA)$ and $L$ is a $1$-cocycle of $A$, then $\phi_{L,a}$ is a morphism of coalgebras. If $\cA$ is a Hopf algebra, then $\phi_{L,a}$ is a morphism of Hopf algebras.
\end{theorem}

\begin{proof}
We will suppress the subscript and write $\phi := \phi_{L,a}$ for ease of notation. The proof is fundamentally the same as the proof of Theorem 2 in \cite{connes1999hopf}, which describes the universality of the pair $(\cH_{CK}, B^+)$.

  We first define $\phi(\bullet) = a$. For a forest $T_1 \cdots T_n$, we have, by multiplicativity,
  \[ \phi(T_1 \cdots T_n) = \phi(T_1) \cdots \phi(T_n). \]
  Any tree with more than one vertex can be written as the result of applying $B^\circ$ to a forest $F$ of smaller trees,
  \[ \phi( B^\circ(T_1 \cdots T_n)) = L( \phi(T_1) \cdots \phi(T_n) ), \]
  so that $\phi \circ B^\circ = L \circ \phi$.
  
  Finally we verify that $\phi$ is a morphism of bialgebras under the additional assumption that $a \in \Prim(\cA)$ and $L$ is a cocycle of $\cA$.  We first verify that this is true for $\bullet \in \cH_B$:
  \[ \Delta_B(\phi(\bullet)) = \Delta_B(a) = a \otimes 1 + 1 \otimes a = (\phi\otimes \phi)(\bullet\otimes1 + 1\otimes\bullet) = (\phi \otimes \phi) \Delta(\bullet)\]
  by the primitivity assumption on $a \in \cA$. Now we will verify it holds for $B^+(T_1 \cdots T_n)$, assuming inductively that it holds for all smaller trees. We have
  \begin{align*}
    \Delta_B(\phi(B^+(T_1 \cdots T_n)))
    &= \Delta_B(L(\phi(T_1 \cdots T_n))) \\
    &= (\id \otimes L) \Delta_B(\phi(T_1 \cdots T_n)) + L(\phi(T_1 \cdots T_n)) \otimes 1 \\
    &= (\id \otimes L) (\phi \otimes \phi) \Delta_{CK}(T_1 \cdots T_n) + L(\phi(T_1 \cdots T_n)) \otimes 1 \\
    &= (\phi \otimes L \circ \phi) \Delta_{CK}(T_1 \cdots T_n) + L(\phi(T_1 \cdots T_n)) \otimes 1 \\
    &= (\phi \otimes \phi \circ B^+) \Delta_{CK}(T_1 \cdots T_n) + \phi(B^+(T_1 \cdots T_n)) \otimes 1 \\
    &= (\phi \otimes \phi) (\id \otimes B^+) \Delta_{CK}(T_1 \cdots T_n) + \phi(B^+(T_1 \cdots T_n)) \otimes 1 \\
    &= (\phi \otimes \phi) \left( (\id \otimes B^+) \Delta_{CK}(T_1 \cdots T_n) + B^+(T_1 \cdots T_n) \otimes 1 \right) \\
    &= (\phi \otimes \phi) \left( \Delta_{CK} B^+(T_1 \cdots T_n) \right)
  \end{align*}
  where the last equality follows from $B^+$ satisfying the $1$-cocycle condition for $(\cH_{CK}, \Delta_{CK})$.
  
  Now consider the coalgebra $\cA$ with coproduct $\Delta'$ and counit $\eps'$. To conclude that $\phi$ is a bialgebra morphism we must show also that
  \[ \eps' \circ \phi = \eps_B, \] 
  where $\eps_B$ is the counit of $\cH_B$, which satisfies $\eps_B(1) = 1$ and $\eps_B(F) = 0$ for nontrivial forests. We recall that $(\eps' \otimes \id_\cA) \circ \Delta' = \id_\cA$ by coalgebra properties and follow the proof given in~\cite{foissy2013introduction}. We also remark that here we are explicitly using the fact that $L$ is a $1$-cocycle of $\cA$.
  
  For any $x \in \cA$,
  \begin{align*}
    L(x) &= (\eps' \otimes \id_\cA)\Delta'(L(x)) \\
         &= (\eps' \otimes \id_\cA)(L(x)\otimes 1 + (\id_\cA \otimes L)\Delta'(x)) \\
         &= \eps'(L(x)) \otimes 1 + (\eps' \otimes L)\Delta'(x) \\
         &= \eps'(L(x)) \otimes 1 + (\id_k \otimes L)(\eps \otimes \id_\cA)\Delta'(x) \\
         &= \eps'(L(x)) \otimes 1 + L(x)
  \end{align*}
  where in the last line we are identifying $x \in \cA$ with $1 \otimes x \in k \otimes \cA$. Thus
  \[ \eps'(L(x)) = 0 \]
  for all $x \in \cA$ and so $\eps' \circ L = 0$.

  Since $a \in \Prim(\cA)$, its counit is automatically zero (see Section 2.3 of \cite{cartier2021classical}). So
  \[ 0 = \eps'(a) = \eps'(\phi(\bullet)) = \eps_B(\bullet) = 0. \]
  If $T$ is any tree with more than one vertex, then $T = B^+(T_1 \cdots T_n)$ for some forest, so we can verify that
  \[ \eps'(\phi(T)) = \eps'(\phi(B^+(T_1 \cdots T_n))) = \eps'(L(\phi(T_1 \cdots T_n))) = 0 = \eps_B(T), \]
  so $\eps' \circ \phi = \eps_B$ for all trees. Finally
  \[ \eps'(\phi(1)) = \eps_B(1) = 1, \]
  since $\eps', \eps_B$ and $\phi$ are all algebra morphisms. Thus we conclude that $\eps' \circ \phi = \eps_B$ and $\phi$ is a morphism of bialgebras.
\end{proof}

\begin{example}
  By choosing $(\cA, L) = (\cH_{red}, B^\circ_{red})$ and $\bullet \in \Prim(\cH_{red})$ as the primitive element, the projection map $\pi$ is induced by the universal property.
\end{example}

\begin{example}
  The first Hochschild cohomology of $k[X]$ is generated by the integration operator, $L(X^n) = \int_0^X t^n \, dt$, where $L(1) = X \neq 0$. Since any $1$-coboundary $A$ satisfies $A(1) = 0$, no $1$-cocycle in $k[X]$ has any representative where the image of $1$ is zero. Thus no map $\phi_{L,a}$ can intertwine $B^\circ$ with the integration operator of $k[X]$.
\end{example}

\section*{Acknowledgements}

I thank my advisor Matilde Marcolli for bringing this topic to my attention, and for her generous support. I am grateful to Lo\"ic Foissy for taking the time to write out folklore results in combinatorial Hopf algebras for my reference and for reading an early draft of this manuscript. No AI was used in the production of this article.

\bibliography{biblio}
\bibliographystyle{plain}

\end{document}